\documentclass{article}
\usepackage{amsfonts}
\usepackage{float}
\usepackage{amsfonts}
\usepackage{amssymb}
\usepackage[table]{xcolor}

        \newtheorem{theorem}{Theorem}

        \newtheorem{remark}[theorem]{Remark}

\usepackage{epsfig}

\title{A Monte Carlo approach to computing stiffness matrices arising 
in polynomial chaos approximations}

\author{Juan Galvis\footnote{Departamento de Matematicas, Universidad Nacional de Colombia, Bogot\'a, Colombia, {\tt jcgalvisa@unal.edu.co}}
 \and O. Andr\'es Cuervo\footnote{Departamento de Matematicas, Universidad Nacional de Colombia, Bogot\'a, Colombia, {\tt oacuervof@unal.edu.co} }}

\begin{document}
\maketitle
\begin{abstract}

We use a Monte Carlo method to assemble finite element matrices for polynomial 
Chaos approximations of elliptic equations with random coefficients. 
In this approach, all required expectations are approximated by  a Monte Carlo method.
The resulting methodology requires dealing with sparse block-diagonal matrices instead of block-full matrices. This leads to the solution of a coupled system of elliptic equations where the coupling is given by a Kronecker  product matrix involving polynomial evaluation matrices.  This generalizes the Classical Monte Carlo approximation and Collocation method for approximating functionals of solutions of these equations. 
\end{abstract}

\section{Introduction}
We consider the computation of the solution of the following random equation,
\begin{equation} \label{problem}
\left\{\begin{array}{rcl}
-\nabla_{x}\cdot (\kappa(x,\cdot)\nabla_{x}
u(x,\cdot) )&=&f(x,\cdot), \mbox{ for all $x\in D$}\\
u(x,\cdot)&=&0, \mbox{ for all  $x\in \partial D$,}
\end{array}\right.
\end{equation}
where $\kappa(x,\cdot) $ is a random field. The forcing $f$  is also allowed to be random.
See ~\cite{Babuska2,Collocation, MK,FP3,Roman,Schwab,Xiu1} and references  
therein. 

Among the classical common  approaches to computing approximations of the solution  
of (\ref{problem}) we can mention Monte Carlo approximations, 
collocation  approaches,  and methods based on 
especial functions expansions such as polynomials chaos methods. 

Briefly speaking, if we know the distribution of the
processes modeling  the coefficient  of the equation, we can
generate some samples of the involved processes and, for each sample,   apply a finite element method (FEM) to obtain an approximation
of the solution for that particular realization. Then, Monte Carlo
approximations   of functionals of the solution 
are  of the form
$\mathbb{E}[g(u(x,\cdot))]\approx\frac{1}{M}
\sum_{i=1}^{M}{ g \big( u^{(a)}}
(x,\omega_{i})\big)$, where $\mathbb{E}~$ denotes the expectation
operator, $M$ is the number of realizations, and
${\hat{u}}^{(a)}(x,\omega_{i})$ is a finite element approximation of
the solution at $x$ for the $ith$-sample $\omega_i$.
This procedure is however very time consuming as it
involves assembling and solving
large linear systems as many times as
trajectories are simulated. 

As an alternative to the Monte Carlo approach, we can use a Chaos expansion 
based method. In this case, we have
 \begin{equation}
 \hat{u}(x,\omega)\approx \hat{u}^{(a)}(x,\omega)=\sum_{\alpha\in 
{\widetilde{\mathcal{I}}}}
 {\hat{u}^{(a)}}_{\alpha}(x)Y_{\alpha}(\omega),
 \label{eq:approximation}
\end{equation}
 where ${\widetilde{\mathcal{I}}}$ is a 
finite index set and  and
$\left\{Y_{\alpha}\right\}_{\alpha\in 
{{\mathcal{I}}}}$ 
is a collection of random variables with known 
probability distributions.
The approximation of functionals  of the solution can be computed directly 
by hand calculations or  we could  use a Monte Carlo
 method based on (\ref{eq:approximation}).
Using this procedure  we need to solve a very large linear system only once. 
This linear system may be hard to deduce and assemble due to the fact that 
is not easy to manage this kind of expansion for many practical cases. The 
effectiveness of this procedure depends mainly on:
 
\begin{enumerate}
\item The kind of expansion used 
in  (\ref{eq:approximation}).
This choice depends on the random variables involved in the definition of the 
random coefficient. If normal random variables are involved, a  \emph{Wiener-Chaos} expansion may be considered. 
\item   The 
finite dimensional 
problem involved in the  computation 
of  the coefficients 
$\{\hat{u}^{(a)}_{\alpha}\}_{\alpha\in{\widetilde{\mathcal{I}}}}$
in  
 (\ref{eq:approximation}).
Usually a Galerkin or Petrov-Galerkin type problem
that uses the original coefficient $\kappa$ or an approximation of it.
\end{enumerate}

One of the main difficulties when implementing polynomial chaos 
finite element methods is that 
stiffness matrix is difficult to obtain and assemble. For instance, for polynomial chaos finite
 element method the expectations appearing  in the bilinear forms can only be computed
  exactly in few special cases - usually for coefficients with explicitly available  expansion  in terms of 
Fourier-Hermite polynomials.
In this paper, we use a Monte Carlo methodology to assemble these matrices.
That is, the expectation appearing in the involved bilinear forms are approximated 
using a Monte Carlo method.  We discuss possible advantages of using this approach 
in comparison with usual approaches. In particular, we find that using Monte Carlo approximation of expectations in bilinear forms might be very advantageous in 
cases where explicit expansions of parameters of the equations are not available and 
especially in the case where these expansions involve functions with small support across 
the domain.  For instance, we mention the case of log-normal coefficients with sparse KL expansions  or with KL expansion  with compact support coefficients.   \\

\section{Problem and proposed approach}
The weak form of (\ref{problem}) is to find $u\in H^1_0(D)\otimes L^2(\mu)$ such that 
\begin{equation}\label{vf}
\int_{D\times \Omega} \kappa \nabla u \nabla v = \int_{D\times \Omega} fv \quad 
\mbox{for all } v \in H^1_0(D)\otimes L^2(\mu).
\end{equation}
The analysis of this variational problem depends on the properties of $\kappa$ 
and $\mu$. We refer the reader to \cite{GalvisSarkisRegularity}  for a review of different approaches to analyzing
this formulation.\\

We assume that the coefficient $\kappa$ and the forcing $f$
are  of the form 
\begin{equation}\label{formofkappa}
\kappa(x,\omega)=\kappa(x,y_1(\omega),y_2(\omega), \dots, y_N(\omega))
\end{equation}
and 
\begin{equation}\label{formoff}
f(x,\omega)=f(x,y_1(\omega),y_2(\omega), \dots, y_N(\omega))
\end{equation}
where $\{ y_j\}_{1\leq j\leq N}$ are radon variables defined in $\Omega$. In this case we can write the solution as $u=u(x,y_1(\omega),y_2(\omega), \dots, y_N(\omega))=u(x,y_1,y_2, \dots, y_N)$. From now one we denote 
$y=(y_1,y_2, \dots, y_N)$.

\begin{remark}
In some practical cases it is given the expansion 
\[\log \kappa(x,\omega)= \sum_{i=1}^N a_i(x) y_i(\omega),
\]
with the functions $a_i$ being of compact support. 
\end{remark}

Introduce the space, 
\begin{equation}\label{discretespace}
\mathcal{P}^{h,M}= P_0^1(\mathcal{T}^h)\otimes P^M(y).
\end{equation}
Here $P_0^1(\mathcal{T}^h)$ is the finite element space of piecewise constant 
functions that vanish on the boundary $\partial D$. The space 
$P^M(y)$ is the space   of polynomials in the variables 
$y_1,y_2, \dots, y_N$ that have a total degree at most $M$.

The discrete problem to approximate (\ref{vf}) is to find 
$u\in \mathcal{P}^{h,M}$ such that 
\begin{equation}\label{gf}
\int_{D\times \Omega} \kappa \nabla u \nabla v = \int_{D\times \Omega} fv \quad 
\mbox{for all } v \in \mathcal{P}^{h,M}.
\end{equation}

Let $\{ \phi_i\}$ be the standard basis of $P_0^1(\mathcal{T}^h)$ and also let 
$\{v_j\}$ be a basis of $P^M(y)$. Consider the set 
$\{ \phi_iv_j\}$ which forms a basis for the space $\mathcal{P}^{h,M}$. We 
reorder the basis functions to $\{ \Phi_I\}$ where $I=(i_1,i_2)$ and 
\[
 \Phi_I(x,y)=\phi_{i_1}(x)v_{i_2}(y).
\]

The matrix form of the problem is written as 
\begin{equation}\label{mf}
\mathcal{A}U=F
\end{equation}
where $\mathcal{A}=(a_{IJ})$ with 
\[
a_{IJ}=\int_{D\times \Omega} \kappa \nabla \Phi_I \nabla \Phi_J=
\int_{D} \left( \int_{\Omega} \kappa  v_{i_2}v_{j_2} \right) \nabla \phi_{i_1} \nabla \phi_{j_1}.
\]
The right hand side vector is given by $F=(f_I)$ where 

\[
f_{I}=\int_{D\times \Omega} f \Phi_I=
\int_D \left( \int_{\Omega} f  v_{i_2}\right)\phi_{i_1}.
\]
Note that in order to assemble the resulting stiffness matrix and load vector  we need to compute expectations with respect to the measure involved. These expectations can be computed using a Monte Carlo method. More precisely,  let us generate samples  
$y^{(1)},y^{(2)}, \dots, y^{(S)}$ of the random vector $y=(y_1,y_2,\dots, y_M)$.
We use the approximation
\begin{equation}\label{mfapprox}
\widetilde{\mathcal{A}}U=F,
\end{equation}
where $\widetilde{\mathcal{A}}=(\widetilde{a}_{IJ})$  and $\widetilde{F}=(\widetilde{f}_I)$ with 
\[
a_{IJ}\approx \widetilde{a}_{IJ}=
\int_{D} \left( \frac{1}{S}\sum_{r=1}^S  \kappa(x,y^{(r)})v_{i_2}(y^{(r)})v_{j_2}(y^{(r)}) \right) \nabla \phi_{i_1} \nabla \phi_{j_1},
\]
and
\[
f_{I}\approx \widetilde{f}_{I}=
\int_{\Omega} \left( \frac{1}{S}\sum_{r=1}^S f(x,y^{(r)})  v_{i_2}(y^{(r)}) \right)\phi_{i_1}.
\]

In this way, the only valuation of coefficient $\kappa$ and forcing $f$ are needed. No especial forms or expansions have to be computed. 

Note that
\[
\widetilde{a}_{IJ}=\frac{1}{S}\sum_{r=1}^S  
\left(\int_{D} \kappa(x,y^{(r)})  \nabla \phi_{i_1} \nabla \phi_{j_1} \right) v_{i_2}(y^{(r)})v_{j_2}(y^{(r)}).
\]
Introduce the matrices
\[
A^{(r)}= (a_{i_1j_1}^{(r)}) \mbox{ where }  a_{i_1j_1}^{(r)}=\int_{D} \kappa(x,y^{(r)})  \nabla \phi_{i_1} \nabla \phi_{j_1},
\]
and 
\[
Z=(z_{ri_2}) \mbox{ where } z_{ri_2}=v_{i_2}(y^{(r)}).
\]
Let $I$ denote the identity matrix of the same size of $B^{(r)}$ and let $V$ be defined 
by the Kronecker product $V=Z\otimes I$. Following \cite{MR2861658}, we have, 
\[
V= \left(\begin{array}{cccc}
v_{1}(y^{(1)})I & v_{2}(y^{(1)})I & \dots & v_{N_2}(y^{(1)})I \\
v_{1}(y^{(2)})I & v_{2}(y^{(2)})I & \dots & v_{N_2}(y^{(2)})I \\
\vdots &\vdots & \ddots & \vdots \\
v_{1}(y^{(S)})I & v_{2}(y^{(S)})I & \dots & v_{N_2}(y^{(S)})I 
\end{array}\right).
\]
Here $N_2=\mbox{dim} P^M$. Introduce also the vectors $F^{(r)}=(f_{i_1}^{(r)})$ where \[f_{i_1}^{(r)}=
\int_D f(x, y^{(r)})\phi_{i_1}.\]

We see that 
\begin{equation}\label{Atilde}
\widetilde{A} = \frac{1}{S} V^T \mbox{diag}(A^{(1)},A^{(2)},\dots, A^{(S)}) V,
\end{equation}
and 
\[
\widetilde{F}=\frac{1}{S} V^T \left(\begin{array}{c} F^{(1)}\\ F^{(2)}\\\vdots\\F^{(S)} \end{array}\right).
\]
Therefore, solving the linear system using an iteration requires only to manage $S$ sparse matrices of usual sparsity pattern and a (possible dense)  $S\times \mbox{dim}(P ^M)$ matrix $Z$. The matrix $Z$ recovers the valuation of the basis of the space $P^M$ at the samples of the 
processes involved. 
Recall that in classical chaos finite element procedure the resulting matrix is blocked dense where each block is a usual finite element matrix.

Note that when $S= \mbox{dim}P^M$ and the basis is selected to be the 
Lagrange polynomials based on the random samples, then, matrix $V$ is the 
identity matrix and we recover the classical Monte Carlo approximation.
Also, when the random samples are collocated in the zeros of suitable orthogonal polynomials, then the resulting method can be viewed as a polynomial chaos method where the expectations are computed using integration rules.

We  note that a similar approach can be implemented starting with the first order formulation 
\begin{equation} \label{problem1order}
\left\{\begin{array}{rcl}
\kappa(x,\cdot)^{-1}q(x,\cdot)&=&-\nabla_{x}
u(x,\cdot) \\
\nabla_{x}\cdot q(x,\cdot)&=&f(x,\cdot), \mbox{ for all $x\in D$}\\
u(x,\cdot)&=&0, \mbox{ for all  $x\in \partial D$.}
\end{array}\right.
\end{equation}
To stress the difference of these two approaches note that when $M=0$, the problem (\ref{mfapprox}) computes the finite element approximation of the solution of the equation with mean coefficient,
\begin{equation} 
\left\{\begin{array}{rcl}
-\nabla_{x}\cdot (\frac{1}{S}\sum_{r=1}^S  \kappa(x,y^{(r)})\nabla_{x}
u(x,\cdot) )&=&f(x,\cdot), \mbox{ for all $x\in D$}\\
u(x,\cdot)&=&0, \mbox{ for all  $x\in \partial D$,}
\end{array}\right.
\end{equation}
while in the case of using the first order system of equations we get an approximation of the solution of the equation  with harmonic mean coefficient, 
\begin{equation} 
\left\{\begin{array}{rcl}
-\nabla_{x}\cdot \left( \left(\frac{1}{S}\sum_{r=1}^S  \kappa(x,y^{(r)})^{-1}\right)^{-1}\nabla_{x}
u(x,\cdot) \right)&=&f(x,\cdot), \mbox{ for all $x\in D$}\\
u(x,\cdot)&=&0, \mbox{ for all  $x\in \partial D$,}
\end{array}\right.
\label{eq:original}
\end{equation}
We refer to \cite{MR2740772} for discussions on these two models.

\section{Some issues of the proposed approach}

We first note that this approach is a generalization of Monte Carlo and Collocation:
\begin{itemize}
\item {\bf Classical Monte Carlo method:} Consider the case where we start with samples $y^{(1)},y^{(2)}, \dots, y^{(S)}$ and we use as a basis of  $P^M(y)$  the Lagrange 
cardinal functions $\{ v_j\}$ based on these samples (that is we have 
$v_j(y^{(k)}) =\delta_{jk}$). In this case we have that $V$ is the identity matrix and therefore $\widetilde{\mathcal{A}}= \mbox{diag}(A^{(1)},A^{(2)},\dots, A^{(S)})$. We recover the Monte Carlo method where $S$ systems are solved independently of each other.
\item {\bf Classical collocation method:}  In case  the basis of  $P^M(y)$ is given by orthogonal polynomials $\{ v_j\}$ and 
 samples $y^{(1)},y^{(2)}, \dots, y^{(S)}$ are collocated at the zeros of these orthogonal polynomials we recover the classical Collocation method. 

\end{itemize}

In the general case of an arbitrary  polynomial base  $\{ v_j\}$ and random samples $y^{(1)},y^{(2)}, \dots, y^{(S)}$ system (\ref{mfapprox})
offers a coupling between otherwise independents solves in the classical Monte Carlo approach. The introduced coupling is generated by a Kronecker product of matrices related to valuation of the polynomial basis at the random samples. \\

Note that the dimension of the polynomial space and the number of samples are, in general, independent of each other. In the case of $S\to  \infty$, then the solution of system (\ref{mfapprox}) will approximate the Galerkin solution obtained with such a
polynomial space for the random variable. The convergence properties when $S\to\infty$ are expected to be those of classical Monte Carlo method. That is, the order of convergence is expecte to be $S^{-1/2}$.

Linear system (\ref{mfapprox}) (where expectations in bilinear forms are approximated using a Monte Carlos approach) offers several interesting issues when compared to 
Linear system (\ref{mf}) (which is the exact bilinear form given by the Galerkin chaos expansion method).
\begin{itemize}
\item {\bf Storage:} In term of storage, the linear system  (\ref{mf}) requires 
to store a $M\times M$ block-dense linear system  where each block corresponds to a sparse matrix of size $O(h^{-d})$ if $x\in R^d$. On the other hand, the linear system
(\ref{mfapprox}) requires the storage of  $S$ sparse matrices of size $O(h^{-d})$ plus a
dense  matrix of size $S\times M$ (or the posibility of performing the evaluations of the fly). 
\item {\bf Local matrices:} In terms of computation of local matrices the assembling of 
the system (\ref{mf}) requires, in general, the computation of $M\times M$ local matrices in each element of the triangulation. Each of these local matrices involves computing 
a coefficient at the quadrature points. For instance, it is usual when using 
polynomial chaos based on Hermite polynomials that  the computation of each coefficient 
may involve sums and products of quantities involving moderate and large numbers such as factorials 
or combinations.  On the other hand system (\ref{mfapprox}) requires, 
in general, $S$ local matrices in each element. The local coefficient of these matrices 
are samples of the coefficient in the original problem so they usually require few function evaluations.

\item {\bf Local matrices for special from of parameters:} Consider the case 
\[ \kappa(x,\omega)= \sum_{i=1}^M a_i(x) y_i(\omega)
\]
with the functions $a_i$ being of compact support (\cite{bachmayr2017sparse}). In this case system  (\ref{mfapprox}) requires  to compute only few local matrices (as it would be for the classical Monte Carlo method). In fact, for an element $\tau\in T^h$, there is only need to compute $\#\{ i : \quad support(a_i)\cap \tau \not = \mbox{empty} \}.$

\item{\bf Error estimates.} The expected bound for the a priori error estimates for the 
solution of (\ref{mf}) using the space  $\mathcal{P}^{h,M}$   (under usual assumptions 
of regularity of forcing term and solution  -\cite{GalvisSarkisinfsup, 
GalvisSarkisRegularity}) is of the form
of $$ || u-u^{h,M}||_{H^1(D)\otimes L^2(\Omega)}\leq C( M^{-1}+h).$$ Using the
 Strang lemma and standard Monte Carlo approximation results, for the solution of 
(\ref{mfapprox}) using $S$ samples to approximate the expectations, one can expect 
and error estimates of the order of $$|| u-u^{h,M,S}||_{H^1(D)\otimes 
L^2(\Omega)}\leq C( M^{-1}+h +S^{-1/2}).$$ We recall that for the standard Monte 
Carlo method for the approximation  of the mean of the solution it is obtained an  error in the form 
\[
\left|\int_{\Omega} u^{h}(x,\omega) - \frac{1}{S} \sum_{r=1}^S 
 u^{h}(x,\omega_i) \right|\leq 
C S^{-1/2}
\]
where $u^h(\cdot,\omega)$ is the solution approximation with mesh size $h$ for the $\omega$ sample.

\item{\bf Iterative method and preconditioning.} Note that, in general, it is hard to construct 
preconditioners for the full block system (\ref{mf}). The form of the final linear system 
(\ref{mfapprox}) with matrix (\ref{Atilde}) (or its saddle point equivalent formulation) is suitable for constructing preconditioners. The construction of robust iterative solvers for 
system (\ref{mfapprox}) is object of future research.
\end{itemize}

\section{Numerical experiments and discussions}
In this  section we present some simple numerical experiments to support the 
ideas proposed before.  As a study case, we approximate the solution of the equation  $$\left\{ \begin{array}{lc}
\mbox{Find } u:[0,1] \rightarrow \mathbb{R} \mbox{ such that }\\
             -(e^{c(x,\omega)}u_x(x,\omega))_x=f(x,\omega) & \mbox{for all} \ \ x\in [0,1],\\
             u(x,\omega)=0 & \mbox{ for } x \in \{ 0,1\},
             \end{array}
   \right.
$$ in two different  ways. For the first case, we consider  the coefficient $$c(x, \omega)=c(x,y)=a(x)y(w)$$ according  to a normal standard  random variable $y$  and $a(x)= \sin(x)$  for $x \in [0,1].$   As the basis of $\mathcal{P}^{h,M}$  we 
$$ v_n(y) = H_n (y) $$  where $H_n$ is the Hermite polynomial of degree $n$. We consider the exact solution,  
$$ u(x, \omega)= u(x, y) = \frac{x(1-x)}{2}e^{-a(x)y}, $$ and therefore
\begin{eqnarray*}
f(x, y)
= H_0(y) + \left(\frac{(1-2x)\mbox{cos}(x)}{2} - \frac{x(1-x)\mbox{sin}(x)}{2} \right)H_1(y).
\end{eqnarray*}
We compare the computed expected value of the solutions with the expected value of the exact solution  that is given by, 
$$ u_0(x)=\displaystyle\int_{\mathbb{R}}u(x,\omega)d\mu(\omega) =  \frac{x(1-x)}{2} e^{\frac{a(x)^2}{2}}.$$  

Denote by $u^{(a)}_{0}(x)$ the approximated mean value. 
We use the $H^1$ error given by  $$ \varepsilon_{H^1} = \displaystyle\sqrt{({u}_{0} -u^{(a)}_{0})^T \cdot A \cdot ({u}_{0} -u^{(a)}_{0})}.$$ 
Note that the error is computed by using the exact solution of the problem. Here the matrix $A$ is defined by  $$ A (i,j) = \int_0^1 \phi'_i(x) \phi'_j(x) dx, $$ with $\phi_i$ and $\phi_j$ the basis functions of finite element method developed on the interval $[0,1]$. \\

Tables \ref{table1} and \ref{table2} show the error behavior as the number of terms of the degree of the polynomials  ($n$) as well as the number of samples ($S$) in the approximation of expectations. Additionally, the last row of the table shows the error  in the approximation of the mean by Monte Carlo method with $S$ samples.  For comparison we only report the error of the computation of the mean (that in  the case of the basis being the Hermite polynomials corresponds to the coefficient of  $H_0$). We observe that the error 
of the computation of the mean of the solution of of the same order as that of the Monte Carlo approximation with some improvement for some higher oder polynomials. \\

\begin{figure}[ht]
\centering
   \includegraphics[width = .48\textwidth, keepaspectratio = true]{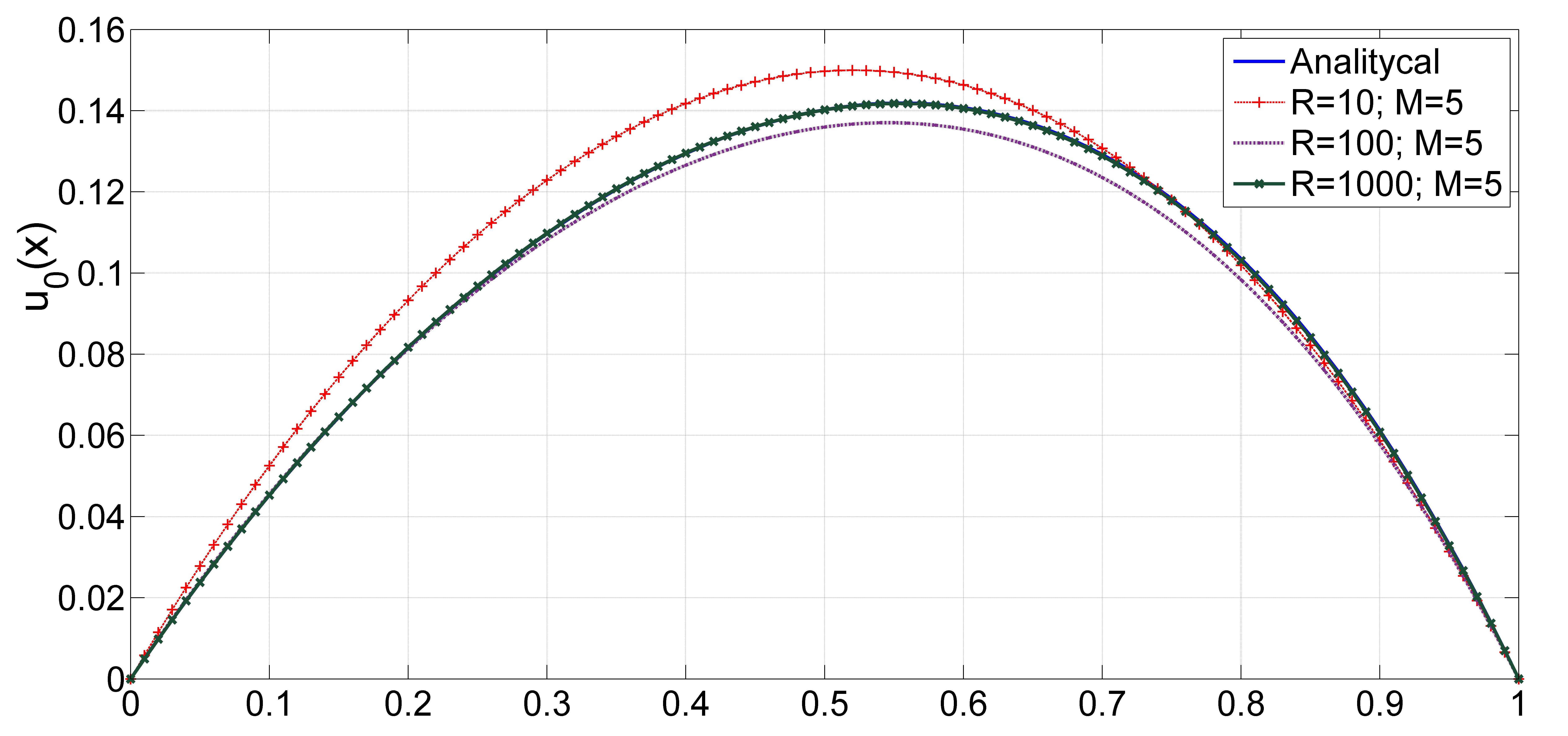}
   \includegraphics[width = .48\textwidth, keepaspectratio = true]{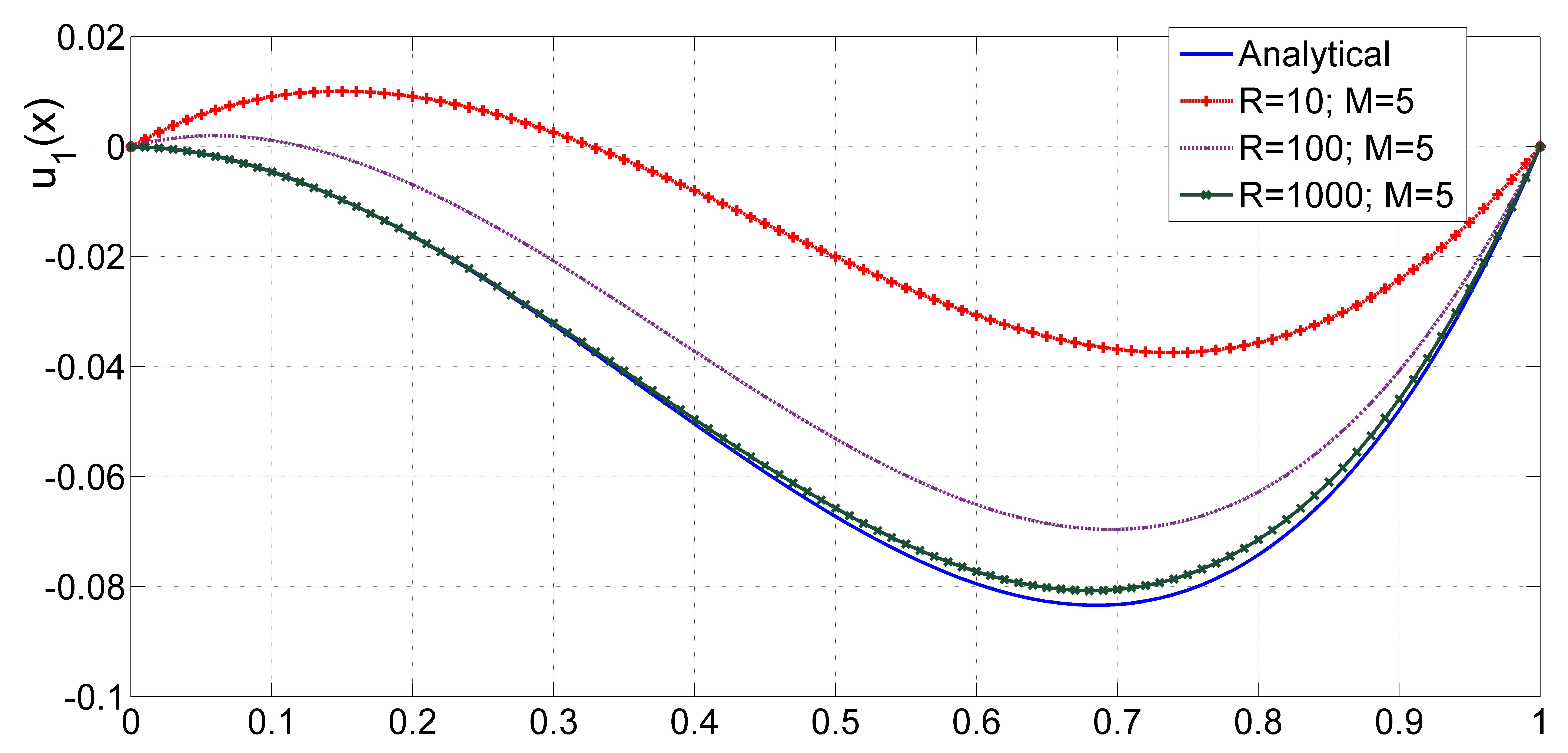}
   \includegraphics[width = .48\textwidth, keepaspectratio = true]{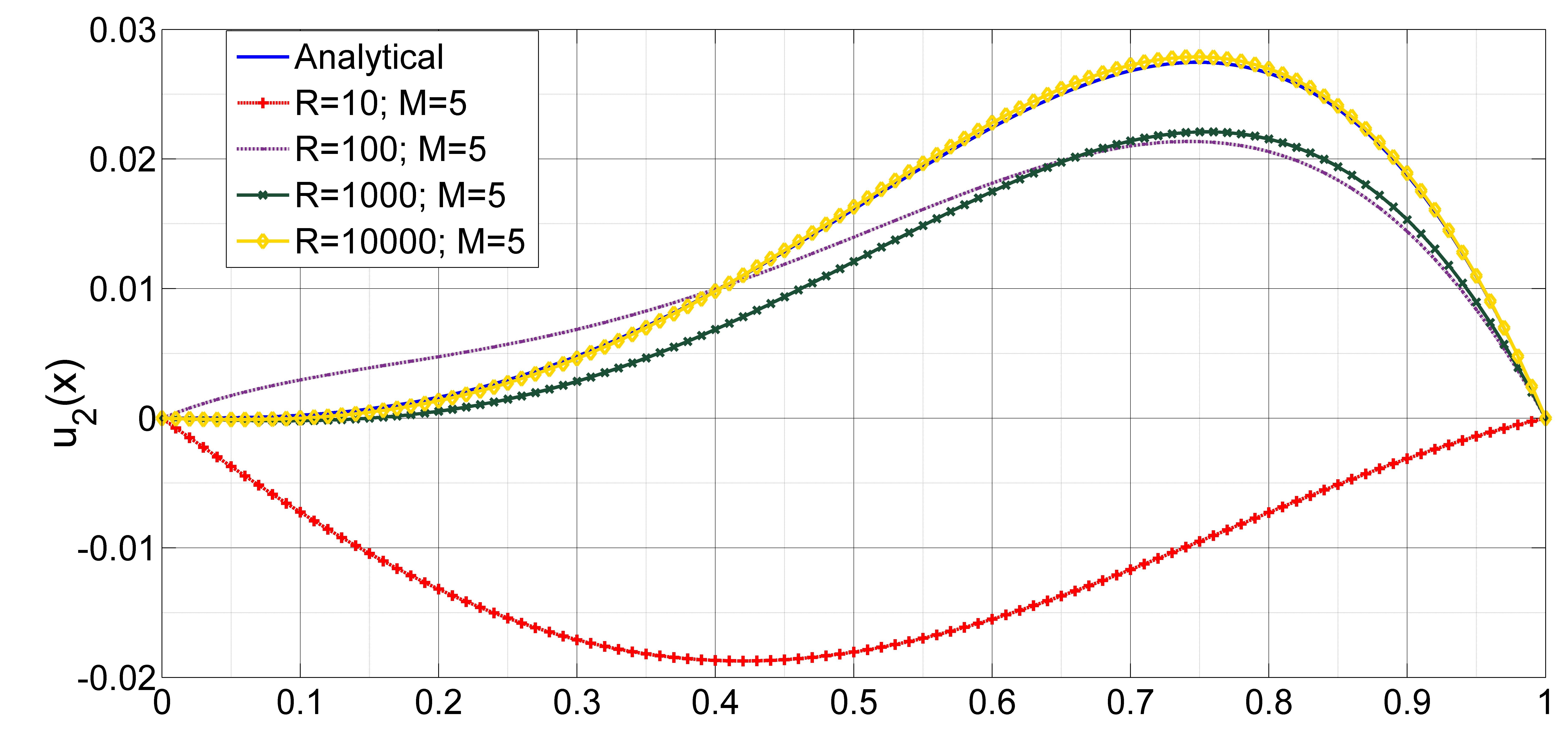}
      \includegraphics[width = .48\textwidth, keepaspectratio = true]{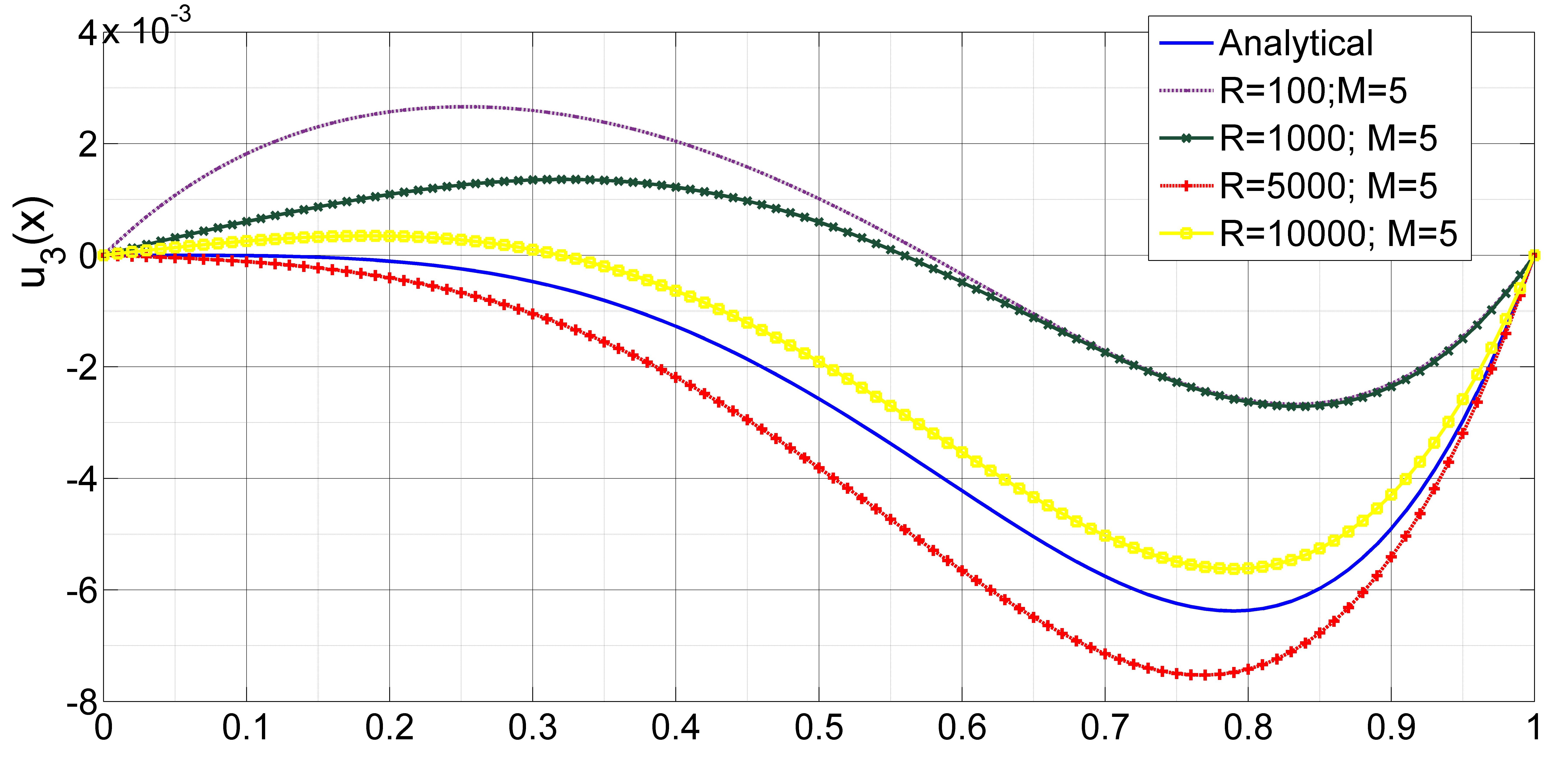}
\caption{Computation of the first four coefficients in the Chaos expansion for a fixed mesh and different number of realizations.}\label{figure1}
\end{figure}

In Figure \ref{figure1}  we illustrate the convergence with respect to the number of samples used in the computation of the stiffness matrix.\\

\begin{table}
\begin{center}
 \begin{tabular}{|c|c|c|c|c|c|c|c|}\hline
\small{\small{$n$}}$\diagdown
$\small{\small{$S$}}&\textbf{100}&\textbf{500}&\textbf{1000}&\textbf{5000} &\textbf{10000} \\ \hline
\textbf{1}  &  0.09721516  &0.10894867& 0.11188709  &  0.11554286  & 0.11545114 \\
\textbf{2}  &  0.03233329  &0.02578792& 0.02564406  &  0.02563924 & 0.02462266 \\
 \textbf{3} & \color{blue}0.01606539  &0.00388315& 0.00398327  &  \color{blue}0.00406651 & \color{blue}0.00421054 \\
 \textbf{4} &\color{blue}0.01493936  &\color{blue}0.00136461& \color{blue}0.00160280  &  \color{blue}0.00559362 & \color{blue}0.00560562 \\
\textbf{5}  &  \color{blue}0.01498601  &\color{blue}0.00127964& \color{blue}0.00155254  &  0.00637033 & 0.00587425 \\ 
 \textbf{6} &  0.01797206  &\color{blue}0.00096630& \color{blue}0.00141063  &  0.00701617 & 0.00601310 \\ \hline
  \textbf{Error MC} &  0.01797206  & 0.0024618    &  0.00266236   &    0.00617800 & 0.00568294 \\ 
 \hline \end{tabular} \\[0.5ex]
  \caption{Error table with Hermite polynomials expansions and $ a (x) = \mbox {sin} (x) $, with error $ \ \varepsilon_ {H ^ 1} $ and $N=100$ elements.}\label{table1}
\end{center}
\end{table}

\begin{table}
\begin{center}
 \begin{tabular}{|c|c|c|c|c|c|c|c|}\hline
\small{\small{$n$}}$\diagdown
$\small{\small{$S$}}&\textbf{100}&\textbf{500}&\textbf{1000}&\textbf{5000} &\textbf{10000}\\ \hline
\textbf{1} &  0.09723178  & 0.10896572 & 0.11190433 & 0.11556026  &  0.11546850 \\
\textbf{2} &  0.03234042 & 0.02579574& 0.02565230 &   0.02564833 & 0.02463159 \\
\textbf{3} & \color{blue}0.01606686 & 0.00388462& 0.00398515  & \color{blue}0.00407083 & \color{blue}0.00421463 \\
\textbf{4} & \color{blue}0.01494001 &\color{blue}0.00136426 & \color{blue}0.00160298  & \color{blue}0.00559687 & \color{blue}0.00560894 \\
\textbf{5} & \color{blue}0.01498675  &\color{blue}0.00127967 & \color{blue}0.00155292   &  0.00637383  &     0.00587771 \\ 
\textbf{6} & \color{blue}0.01351294  & \color{blue}0.00096802& \color{blue}0.00141142    & 0.00701985  &            0.00601661  \\ \hline
\textbf{Error MC} &   0.01797746  & 0.00246442    &  0.00266460 &  0.00617649     &    0.00568098  \\ \hline \end{tabular} \\[0.5ex]
\caption{Error table with Hermite polynomials expansion and $ a (x) = \mbox {sin} (x) $, with error $ \ \varepsilon_ {H ^ 1} $ and $N=1000$ elements.}\label{table2}
\end{center}
\end{table}

We also consider the $L^2$ error given by, 
$$\varepsilon_{L^2} = \displaystyle\sqrt{(\overline{u}_{0} -u^{(a)}_{0})^T \cdot M \cdot (\overline{u}_{0} -u^{(a)}_{0})} $$
where the matrix $M$ is given by
$$ M(i,j) = \int_0^1 \phi_i(x) \phi_j(x) dx.$$\\
The results are shown in Table \ref{table2}.\\
\begin{table}
\begin{center}
 \begin{tabular}{|c|c|c|c|c|c|c|c|}\hline
\small{\small{$n$}}$\diagdown
$\small{\small{$S$}}&\textbf{100}&\textbf{500}&\textbf{1000}&\textbf{5000} &\textbf{10000} \\ \hline
\textbf{1} &             0.02334552 &            0.02706124&             0.02794990 &             0.02910265 & 0.02908693 \\ 
 \textbf{2}&             0.00717723 &            0.00515693&             0.00498519 &             0.00466146 & 0.00441017\\ 
 \textbf{3}& \color{blue}0.00384046 &            0.00066324&             0.00061879 & \color{blue}0.00094036 & \color{blue}0.00108872\\ 
 \textbf{4}& \color{blue}0.00356213 &\color{blue}0.00024458& \color{blue}0.00025382 & \color{blue}0.00155189 & 0.00157809 \\ 
 \textbf{5}& \color{blue}0.00355041 &\color{blue}0.00020888& \color{blue}0.00023186 &              0.00170032 & 0.00162147 \\
 \textbf{6}& \color{blue}0.00273260 &\color{blue}0.00014588& \color{blue}0.00019812 &              0.00185357 & 0.00165381 \\ \hline
   \textbf{Error MC} & 0.00469386  & 0.00052470 & 0.00050020 & 0.00160435 & 0.00152867\\ \hline
  \end{tabular} \\[0.5ex]
    \caption{Error table with  Hermite polynomials expansion and $ a (x) = \mbox {sin} (x) $, with error $ \ \varepsilon_ {L^2} $ and $N=100$ elements.}\label{table3}
\end{center}
\end{table} 

\begin{table}
\begin{center}
 \begin{tabular}{|c|c|c|c|c|c|c|c|}\hline
\small{\small{$n$}}$\diagdown
$\small{\small{$S$}}&\textbf{100}&\textbf{500}&\textbf{1000}&\textbf{5000} &\textbf{10000} \\ \hline
\textbf{1} &             0.02334907 &             0.02706503 &             0.02795377 &             0.02910659 & 0.02909087 \\ 
 \textbf{2}&             0.00717794 &             0.00515767 &             0.00498599 &             0.00466245 & 0.00441112 \\ 
\textbf{3} & \color{blue}0.00384010 &             0.00066271 &             0.00061843 & \color{blue}0.00094151 & \color{blue}0.00108988 \\ 
 \textbf{4}& \color{blue}0.00356166 & \color{blue}0.00024396 & \color{blue}0.00025340 & \color{blue}0.00155316 & 0.00157937 \\ 
 \textbf{5}& \color{blue}0.00354995 & \color{blue}0.00020843 & \color{blue}0.00023160 &             0.00170164 & 0.00162277 \\
 \textbf{6}& \color{blue}0.00273234 & \color{blue}0.00014648 & \color{blue}0.00019812 &             0.00185493 & 0.00165512 \\ \hline
   \textbf{Error MC} & 0.00469562  & 0.00052583& 0.00050122 & 0.00160350 & 0.00152775 \\ \hline
  \end{tabular} \\[0.5ex]
 \caption{Error table with  Hermite polynomials expansions and $ a (x) = \mbox {sin} (x) $, with error $ \ \varepsilon_ {L^2} $ and $N=1000$ elements.}\label{table4}
\end{center}
\end{table}

We now consider the random coefficient, given by 
$$c(x,\omega)=c(x,y_1,y_2)=a_1(x)y_1(\omega)+a_2(x)y_2(\omega)$$ with $y_1, y_2$ independent  random variables   with normal standard distribution. Also $a_1(x)=\mbox{sin}(x)$ and  $a_2(x)=\mbox{cos}(x)$ for   $x \in [0,1].$ Again, we start from the exact solution
\begin{eqnarray*}
&&u_x(x,y_1,y_2)  =  \left[\frac{1-2x}{2} - \frac{x(1-x)a_1'(x)}{2}y_1 - \frac{x(1-x)a_2'(x)}{2}y_2 \right]e^{-c(x,y_1,y_2) }\\
& = & \left[\frac{1-2x}{2} - \frac{x(1-x)a_1'(x)}{2}y_1 - \frac{x(1-x)a_2'(x)}{2}y_2 \right]e^{-(a_1(x)y_1(w)+a_2(x)y_2(\omega)) },
\end{eqnarray*}   
and 
$$  f(x, y_1 , y_2 ) =  1 +  \left(\frac{(1-2x)\mbox{cos}(x)}{2} - \frac{x(1-x)\mbox{sin}(x)}{2} \right) y_1 $$
\begin{flushright}
$   - \displaystyle\left(\frac{(1-2x)\mbox{sin}(x)}{2} + \frac{x(1-x)\mbox{cos}(x)}{2} \right) y_2.  $ 
\end{flushright}

Considering the random vector  $y=(y_1,y_2)$ and the  collection of random variables $ \{ y_i \}_{i \in \mathcal{I}}$ with $i=(i_1,i_2)$ multi-index, we define $$y_i = H_{i_1}(y_1) H_{i_2}(y_2)$$ with $H_{i_1}, H_{i_2}$ the Hermite polynomials of degrees $i_1$ and $i_2$, respectively.

By direct calculations we have, 
$$ u_0(x)=   \frac{x(1-x)}{2} e^{\frac{a_1(x)^2 + a_2(x)^2}{2}}, $$
as the analytic mean. The resulting error are summarized in the following tables, where now $S$ is the amount of realizations of random vector $y=(y_1,y_2),$ 
  
\begin{table}
\begin{center}
 \begin{tabular}{|c|c|c|c|c|c|c|c|}\hline
$n \setminus S$ &\textbf{100}&\textbf{500}&\textbf{1000} &\textbf{5000}& \textbf{10000} \\ \hline
\textbf{0} & 0.28898445 & 0.30420350 &0.30251536 & 0.30017828 &0.30291653 \\
 \textbf{1}& 0.10133617 & 0.11723215 &0.11672713 & 0.12277728 &0.12391400\\ 
\textbf{2} & \color{blue}0.02882608 & \color{blue}0.01780323 &0.03874657 & 0.03399488 &0.03459619\\
\textbf{3} & 0.06055101 & \color{blue}0.00736567 &0.02053047 &0.00714912 &0.00813146\\
\textbf{4} & 0.07430252 &\color{blue}0.01670585 &0.00974796 & \color{blue}0.00361374 &0.00347825\\
\textbf{5} & 1.62133451 &0.02330246 &\color{blue}0.00243587 & \color{blue}0.00612991 &0.00305979 \\
\textbf{6} & 2.14769303 & 0.02848221&\color{blue}0.00413190 & 0.00721317 &0.00347518 \\
\hline
 \textbf{Error MC} &  0.03446099 &0.01998285 & 0.00903498 &0.00624899 &0.00272246  \\ \hline
  \end{tabular} \\[0.5ex]
  
  \medskip
  
 \begin{tabular}{|c|c|c|c|c|c|c|c|}\hline
$n \setminus S$ &\textbf{100}&\textbf{500}&\textbf{1000} & \textbf{5000} & \textbf{10000} \\ \hline
\textbf{0} & 0.09138038 & 0.09619608 & 0.09566152 &0.09491099 &0.09579028 \\
 \textbf{1}& 0.03180108 &  0.03702665 & 0.03685698 &  0.03874437 &0.03914721\\
\textbf{2} &\color{blue} 0.00469637 & \color{blue}0.00556492 & 0.01205270 & 0.01057680 &0.01080931 \\
\textbf{3} & 0.01910321 & \color{blue}0.00120043 & 0.00646981 & 0.00173947  &0.00216430 \\
\textbf{4} & 0.09078174 & \color{blue}0.00428636 & 0.00308585 & \color{blue}0.00056435 &0.00052493 \\
\textbf{5} & 0.27747279 & 0.00750032 &\color{blue}0.00062702 & \color{blue}0.00138300 &0.00047162 \\
\textbf{6} & 0.88405006 & 0.00904359 &\color{blue}0.00131699 & 0.00185361  &0.00079010  \\
\hline
 \textbf{Error MC} &  0.01093222  & 0.00632597 & 0.00285162 & 0.00151079 &0.00043792  \\ \hline
  \end{tabular} \\[0.5ex]
  \caption{Error table for multiplication of the Hermite polynomials  base, $a_1(x)=\mbox{sin}(x)$, $a_2(x)=\mbox{cos}(x)$, with the error $\varepsilon_{H^1}$ and $\varepsilon_{L^2}$, respectively and $N=100$ elements.}\label{Tabla}
\end{center}
\end{table} 

\section{Conclusions}

We studied the use of a Monte Carlo method to assemble finite element matrices for polynomial Chaos approximations of elliptic equations with random coefficients. 
In this approach, all required expectations are approximated by  a Monte Carlo method.
This leads to the solution of a coupled system of elliptic equations where the coupling is given by a Kronecker  product matrix involving polynomial evaluation matrices.  This generalizes the Classical Monte Carlo approximation and Collocation method for approximating functionals of the solution of these equations. 
The resulting methodology requires dealing with sparse block-diagonal matrices instead of block-full matrices.

\bibliographystyle{alpha}

\def\cprime{$'$}

\end{document}